\documentclass[10pt]{article}
\usepackage{epsf}
\usepackage{epsfig}
\usepackage{color}

\def\PP{\rm \hbox{I\kern-.2em\hbox{P}}}
\def\RR{\rm \hbox{I\kern-.2em\hbox{R}}}
\def\NN{\rm \hbox{I\kern-.2em\hbox{N}}}
\def\ZZ{{\rm {{\rm Z}\kern-.28em{\rm Z}}}}
\def\CC{\rm \hbox{C\kern -.5em {\raise .32ex \hbox{$\scriptscriptstyle
|$}}\kern-.22em{\raise .6ex \hbox{$\scriptscriptstyle |$}}\kern .4em}}

\def\<{\langle}
\def\>{\rangle}

\def\cT{{\cal T}}

\def\cR{{\cal R}}

\def\cQ{{\cal Q}}

\def\cP{{\cal P}}

\def\Chi{\raise .3ex
\hbox{\large $\chi$}} 
\def\lsima{\hbox{\kern -.6em\raisebox{-1ex}{$~\stackrel{\textstyle<}{\sim}~$}}\kern -.4em}
\def\lsim{\hbox{\kern -.2em\raisebox{-1ex}{$~\stackrel{\textstyle<}{\sim}~$}}\kern -.2em}
\def\({\Bigl (}
\def\){\Bigr )}
\def\({\Bigl (}
\def\){\Bigr )}

\newcommand{\be}{\begin{equation}}
\newcommand{\ee}{\end{equation}}
\newcommand{\bea}{$$ \begin{array}{lll}}
\newcommand{\eea}{\end{array} $$}
\newcommand{\bi}{\begin{itemize}}
\newcommand{\ei}{\end{itemize}}
\newcommand{\iref}[1]{(\ref{#1})}
\newtheorem{theorem}{Theorem}[section]
\newtheorem{remark}[theorem]{Remark}
\newtheorem{lemma}[theorem]{Lemma}

\newtheorem{definition}[theorem]{Definition}

\newtheorem{prop}[theorem]{Proposition}

\def\I{{\rm \hbox{I\kern-.2em\hbox{I}}}}
\def\P{{\rm \hbox{I\kern-.2em\hbox{P}}}}
\def\H{{\rm \hbox{I\kern-.2em\hbox{H}}}}
\def\R{{\rm \hbox{I\kern-.2em\hbox{R}}}}
\def\N{{\rm \hbox{I\kern-.2em\hbox{N}}}}
\def\Z{{\rm {{\rm Z}\kern-.28em{\rm Z}}}}
\def\C{{\rm \hbox{C\kern -.5em {\raise .32ex \hbox{$\scriptscriptstyle
|$}}\kern-.22em{\raise .6ex \hbox{$\scriptscriptstyle |$}}\kern .4em}}}

\def\mI{{\mathbb I}}

\def\proof{{\noindent \bf Proof: }}

\def\ve{\varepsilon}
\def\sq{\hfill $\diamond$\\}

\def\R{\mathbb R}

\def\card{\#}

\def\u{{\bf u}}
\def\v{{\bf v}}

\def\sep{ : \;}

\def\seqR{(\cR_N)_{N\geq 0}}

\def\bs{\setminus}

\usepackage{amsmath}
\usepackage{amsfonts}
\usepackage{amssymb}
\usepackage{url}

\DeclareMathOperator\diam{diam}
\DeclareMathOperator\Vect{Vect}
\DeclareMathOperator\interp{{\rm I}}

\definecolor{VioletJB}{rgb}{0.55,0.2,0.90}

\usepackage{stmaryrd}

\def\balpha{ {\bf \alpha} }

\begin{document}
\title{Sharp asymptotics of the $L_p$ approximation error for interpolation on block partitions}
\author{Yuliya Babenko, Tatyana Leskevich, Jean-Marie Mirebeau}
\maketitle
\date{}
\begin{abstract}
{

Adaptive approximation (or interpolation) takes into account local variations in the behavior of the given function, adjusts the approximant depending on it, and hence yields the smaller error of approximation. The question of constructing optimal approximating spline {\it {for each function}} proved to be very hard. In fact, no polynomial time algorithm of adaptive spline approximation can be designed and no exact formula for the optimal error of approximation can be given. Therefore, the next natural question would be to study the asymptotic behavior of the error and construct asymptotically optimal sequences of partitions.

In this paper we provide sharp asymptotic estimates for the error of interpolation by splines on block partitions in $\RR^d$. We  consider various projection operators to define the interpolant and provide the analysis of the exact constant in the asymptotics as well as its explicit form in certain cases.
}
\end{abstract}

\newpage

\noindent
Yuliya Babenko\\
Department of Mathematics and Statistics\\
Sam Houston State University\\
Box 2206\\
Huntsville, TX, USA 77340-2206\\
Phone: 936.294.4884\\
Fax: 936.294.1882\\
Email: babenko@shsu.edu\\

\noindent
Tatyana Leskevich\\
Department of Mathematical Analysis \\
Dnepropetrovsk National University \\
pr. Gagarina, 72, \\
Dnepropetrovsk, UKRAINE, 49050 \\
Email: tleskevich@gmail.com \\

\noindent
Jean-Marie Mirebeau\\
UPMC Univ Paris 06, UMR 7598, Laboratoire Jacques-Louis Lions, F-75005, Paris, France\\
CNRS, UMR 7598, Laboratoire Jacques-Louis Lions, F-75005, Paris, France\\
Email: mirebeau@ann.jussieu.fr\\

\newpage

\section {Introduction}

The goal of this paper is to study the adaptive approximation by interpolating splines defined over block partitions in $\RR^d$.
With the help of introduced projection operator we shall handle the general case, and then apply the obtained estimates to several different interpolating schemes most commonly used in practice.

Our approach is to introduce the ``error function'' which reflects the interaction of approximation procedure with polynomials. Throughout the paper we shall study the asymptotic behavior of the approximation error and, whenever possible, the explicit form of the error function which plays a major role in finding the constants in the formulae for exact asymptotics.



%


\subsection{The projection operator}

Let us first introduce the definitions that will be necessary to state the main problem and the results of this paper.

We consider a fixed integer $d\geq 1$ and we denote by $x=(x_1, \cdots , x_d)$ the elements of $\R^d$. A block $R$ is a subset of $\R^d$ of the form
$$
R = \prod_{1\leq i\leq d} [a_i,b_i] 
$$
where $a_i< b_i$, for all $1\leq i\leq d$. 
For any block $R\subset \R^d$, by $L_{p}(R)$, $1\le p\le \infty$, we denote the space of measurable functions $f:R \to\RR$ for which the value
$$
  \|f\|_p = \|f\|_{L_p(R)} : = \left\{\begin{array}{ll}
                \left(\displaystyle\int\limits_{R} |f (x) | ^p dx\right)^{\frac 1p},   &{\rm if}\;\;\;    1\leq  p < \infty, \\ [10pt]
                {\rm ess sup} \{|f (x) | \sep x \in R \},                                    &{\rm if}\;\;\;       p =\infty.
\end{array}\right.
$$
is finite. 
We also consider the space $C^0(R)$ of continuous functions on $R$ equipped
with the uniform norm $\|\cdot\|_{L_\infty(R)}$.
We shall make a frequent use of the canonical block $\mI^d$, where $\mI$ is the interval
$$
\mI := \left[-\frac 1 2, \frac 1 2\right].
$$
Next we define the space $V := C^0(\mI^d)$ and the norm $\|\cdot\|_V := \|\cdot\|_{L_\infty(\mI^d)}$.
Throughout this paper we consider a linear and bounded (hence, continuous) operator
$
\interp : V\to V. 
$
This implies that there exists a constant $C_{\interp}$ such that
\be
\label{contI}
\|\interp u\|_V \leq C_{\interp} \|u\|_V \text{ for all } u\in V.
\ee
We assume furthermore that $\interp$ is a projector, which means that it satisfies
\be
\label{projAxiom}
\interp\circ \interp = \interp.
\ee
Let $R$ be an arbitrary block. 
It is easy to show that there exists a unique $x_0\in \R^d$ and a unique diagonal matrix $D$ with positive diagonal coefficients such that the transformation
\be
\label{defPhi}
\phi(x) := x_0+ Dx \ \text{ satisfies } \ \phi(\mI^d) = R.
\ee
The volume of $R$, denoted by $|R|$, is equal to $\det(D)$. For any function $f\in C^0(R)$ we then define
\be
\label{interpR}
\interp_R f := \interp(f\circ\phi)\circ \phi^{-1}.
\ee
Note that
\be
\label{changeRect}
\|f-\interp_R f\|_{L_p(R)} = (\det D)^{\frac 1 p}\|f \circ \phi - \interp(f\circ\phi)\|_{L_p(\mI^d)}.
\ee
A block partition $\cR$ of a block $R_0$ is a finite collection of blocks such that their union covers $R_0$ and which pairwise intersections have zero Lebesgue measure.
If $\cR$ is a block partition of a block $R_0$ and if $f\in C^0(R_0)$, by
$
\interp_\cR f \in L_\infty(R_0)
$
we denote the (possibly discontinuous) function which coincides with $\interp_R f$ on the interior of each block $R\in \cR$.

{\bf Main Question.}
The purpose of this paper is to understand the asymptotic behavior of the quantity
$$
\|f-\interp_{\cR_N} f \|_{L_p(R_0)}
$$
{\it {for each given function $f$}} on $R_0$ from some class of smoothness,
where $(\cR_N)_{N\geq 1}$ is a sequence of block partitions of $R_0$ that are optimally adapted to $f$. 

Note that the exact value of this error can be explicitly computed only in trivial cases. Therefore, the natural question is to  study the asymptotic behavior of the error function, i.e. the behavior of the error as the number of elements of the partition ${\cR_N}$ tends to infinity.

Most of our results hold with only assumptions \iref{contI} of continuity of the operator $\interp$, the projection axiom \iref{projAxiom}, and the definition of $\interp_R$ given by \iref{interpR}. Our analysis therefore applies to various projection operators $\interp$, such as the $L_2$-orthogonal projection on a space of polynomials, or spline interpolating schemes described in \S\ref{exampleI}.

\subsection{History}

The main problem formulated above is interesting for functions
of arbitrary smoothness as well as for various classes of splines
(for instance, for splines of higher order, interpolating splines,
best approximating splines, etc.). In the univariate case general
questions of this type have been investigated by many authors. The
results are more or less complete and have numerous applications
(see, for example, ~\cite{LSh}).

Fewer results are known in the multivariate case. Most of them are
for the case of approximation by splines on triangulations (for
review of existing results see, for instance ~\cite{Gr3, KL, us, Cohen, JM}).
However, in applications where preferred directions exist, box
partitions are sometimes more convenient and efficient.

The first result on the error of interpolation on
rectangular partitions by bivariate splines linear in each variable
(or bilinear) is due to D'Azevedo ~\cite{Daz1} who obtained
local (on a single rectangle) error estimates. In ~\cite{PhD} Babenko obtained the exact
asymptotics for the error (in $L_1$, $L_2$, and $L_{\infty}$ norms)
of interpolation of $C^2(\mI^d)$ functions by bilinear splines.

In ~\cite{JAT} Babenko generalized the result to interpolation and
quasiinterpolation of a function $f\in C^2(\mI^d)$ with arbitrary but fixed throughout the domain signature (number of positive and negative second-order partial derivatives). However, the
norm used to measure the error of approximation was uniform.

In this paper we use a different, more abstract, approach which allows us to obtain the exact asymptotics of the error  in a more general framework which can be applied to many particular interpolation schemes by an appropriate choice of the interpolation operator. In general, the constant in the asymptotics is implicit. However, imposing additional assumptions on the interpolation operator allows us to compute the constant explicitly.

The paper is organized as follows. Section \ref{subsecApprox}
contains the statements of main approximation results. The closer
study of the error function, as well as its explicit formulas under
some restrictions, can be found in Section \ref{studyK}. The proofs
of the theorems about asymptotic behavior of the error are contained
in Section \ref{proofTh}.

\subsection{Polynomials and the error function}

In order to obtain the asymptotic error estimates we need to study the interaction of the projection operator $\interp$ with polynomials.

The notation ${\balpha}$ always refers to a  $d$-vector
of non-negative integers
$$
\balpha = (\balpha_1,\cdots , \balpha_d) \in \ZZ_+^d.
$$
For each $\balpha$ we define
the following quantities
$$
|\balpha|:= \sum_{1\leq i\leq d} \balpha_i, \quad \balpha! := \prod_{1\leq i\leq d} \balpha_i!, \quad \max(\balpha) := \max_{1\leq i\leq d} \balpha_i.
$$
We also define the monomial
$$
X^\balpha := \prod_{1\leq i\leq d} X_i^{\balpha_i},
$$
where the variable is $X=(X_1,...,X_d)\in \RR^d$.
Finally, for each integer $k\geq 0$  we define the following three vector spaces of polynomials
\be
\label{defPk}
\begin{array}{lcl}
\P_k    & :=& \Vect\{X^\balpha \sep |\balpha| \leq k\}, \\
\P_k^* &:=& \Vect\{X^\balpha \sep \max(\balpha) \leq k \text{ and } |\balpha| \leq k+1\}, \\
\P_k^{**} &:=& \Vect\{X^\balpha \sep \max(\balpha) \leq k\}.
\end{array}
\ee 
Note that clearly $\dim(\P_k^{**}) = (k+1)^d$.
In addition, a classical combinatorial argument shows that 
$$
\dim \P_k = \binom {k+d}{d}.
$$
Furthermore,
$$
\dim\P_k^* = \dim\P_{k+1} - d = \binom {k+d+1}{d} - d.
$$

By $V_{\interp}$ we denote the image of $\interp$, which is a subspace of $V = C^0(\mI^d)$. Since
$\interp$ is a projector \iref{projAxiom}, we have \be \label{defVI}
V_{\interp} = \{\interp(f): \;\; f\in V\} = \{f \in V : \;\;
f=\interp(f)\}. \ee From this point on, the integer $k$ is
fixed and defined as follows \be \label{defk} k = k(\interp) := \max
\{k'\geq 0 \sep \P_{k'} \subset V_{\interp} \} \ee Hence, the
operator $\interp$ reproduces polynomials of total degree less or
equal than $k$. (If $k=\infty$ then we obtain, using the density of
polynomials in $V$ and the continuity of $\interp$, that $\interp(f)
= f$ for all $f\in V$. We exclude this case from now on.)

In what follows, by $m$ we denote the integer defined by
\be
\label{defm}
m = m(\interp) :=k+1,
\ee
where $k = k(\interp)$ is defined in \iref{defk}.
By $\H_m$ we denote the space of homogeneous polynomials of degree $m$
$$
\H_m := \Vect\{ X^\balpha\sep |\balpha| = m\}.
$$
We now introduce a function $K_I$ on $\H_m$, further referred to as the ``error function''.
\begin{definition}[Error Function]
For all $\pi \in \H_m$ \be \label{defK} K_I(\pi) := \inf_{|R|=1}
\|\pi - \interp_R \pi\|_{L_p(R)}, \ee where the infimum is taken
over all blocks $R$ of unit $d$-dimensional volume.
\end{definition}

The error function $K$ plays a major role in our asymptotical error estimates developed in the next subsection.
Hence, we dedicate \S2 to its close study, and we provide its explicit form in various cases.

The optimization \iref{defK} among blocks can be rephrased into an  optimization among diagonal matrices.
Indeed, if $|R|=1$, then there exists a unique $x_0\in \R^d$ and a unique diagonal matrix with positive coefficients such that $R = \phi(\mI^d)$ with $\phi(x) = x_0+ Dx$. Furthermore, the homogeneous component of degree $m$ is the same in both $\pi\circ \phi$ and $\pi\circ D$, hence $\pi\circ \phi - \pi\circ D \in \P_k$ (recal that $m = k+1$) and therefore this polynomial is reproduced by the projection operator $\interp$. Using the linearity of $\interp$, we obtain
$$
\pi \circ \phi - \interp (\pi \circ \phi) = \pi \circ D - \interp (\pi \circ D).
$$
Combining this with \iref{changeRect}, we obtain that
\be
\label{KD}
K_I(\pi) = \inf_{\substack{\det D = 1\\D\geq 0}} \|\pi\circ D - \interp(\pi\circ D)\|_{L_p(\mI^d)},
\ee
where the infimum is taken over the set of diagonal matrices with non-negative entries and unit determinant.

\subsection{Examples of projection operators}
\label{exampleI}

In this section we define several possible choices for the projection operator $\rm I$ which are consistent with \iref{defk} and, in our opinion, are most useful for practical purposes.
However, many other possibilities could be considered.
\begin{definition}[$L_2(\mI^d)$ orthogonal projection]
\label{L2Proj}
We may define $\interp(f)$ as the $L_2(\mI^d)$ orthogonal projection of $f$ onto one of the spaces of polynomials $\P_k$,
$\P_k^*$ or $\P_k^{**}$ defined in \iref{defPk}.
\end{definition}

If the projection operator $\interp$ is chosen as in Definition \ref{L2Proj}, then a simple change of variables shows that for any block $R$, the operator $\interp_R$ defined by \iref{interpR} is the $L_2(R)$ orthogonal projection onto the same space of polynomials.



To introduce several possible interpolation schemes for which we obtain the estimates using our approach,
we consider a set $U_k\subset \mI$ of cardinality $\#(U_k) = k+1$ (special cases are given below).
For any $\u = (u_1, \cdots u_d) \in U_k^d$ we define an element of $\P_k^{**}$ as follows
$$
\mu_\u (X):= \prod_{1 \leq i\leq d} \left( \prod_{\substack{v\in U_k\\ v\neq u_i}} \frac{X_i - v}{u_i-v}\right) \in \P_k^{**}.
$$

Clearly, 
$
\mu_\u (\u) = \mu_\u (u_1, \cdots , u_d) = 1
$
and 
$
\mu_\u(\v) = \mu_\u(v_1, \cdots , v_d) = 0
$
if $\v = (v_1, \cdots , v_d)\in U_k^d$ and $\v\neq \u$. 

It follows that the elements of  $B := (\mu_\u)_{\u\in U_k^d}$ are linearly independent. Since $\# (B) = \#(U_k^d) = (k+1)^d = \dim(\P_k^{**})$, $B$ is a basis of $\P_k^{**}$.

Therefore, any element of $\mu \in \P_k^{**}$ can be written in the form
$$
\mu(X) = \sum_{\u\in U_k^d} \lambda_\u \mu_\u(X).
$$
It follows that there is a unique element of $\mu\in \P_k^{**}$ such that $\mu(\u) = f(\u)$ for all $\u\in U_k^d$. We define $\interp f := \mu$, namely 
$$
(\interp f)(X) := \sum_{\u\in U_k^d} f(\u) \mu_\u(X) \in \P_k^{**}.
$$
We may take $U_k$ to be the set of $k+1$ equi-spaced points on $\mI$
\be
\label{interpEqui}
U_k = \left\{-\frac 1 2 + \frac n k \sep 0 \leq n \leq k\right\}.
\ee
We obtain a different, but equally relevant, operator $\interp$ by choosing $U_k$ to be the set of Tchebychev points on $\mI$
\be
\label{interpTche}
U_k = \left\{\frac 1 2 \cos\left( \frac {n\pi} k\right) \sep 0 \leq n \leq k\right\}.
\ee
Different interpolation procedures can be used to construct $\interp$. Another convenient interpolation scheme is to take
$$
\interp(f) \in \P_k^*
$$
and $\interp(f) = f$ on a subset of $U_k^d$. This subset  contains $\dim \P_k^*$ points, which are convenient to choose first on the boundary of $\mI^d$ and then (if needed) at some interior lattice points.
Note that since $\dim \P_k^*< \#(U_k^d) =(k+1)^d$, it is always possible to construct such an operator.

If the projection operator $\interp$ is chosen as described above, then for any block $R$ and any $f\in C^0(R)$, $\interp_R(f)$ is the unique element of respective space of polynomials which coincides with $f$ at the image $\phi(p)$ of the points $p$ mentioned in the definition of $\interp$, by the transformation $\phi$ described in \iref{defPhi}.

\subsection{Main results}
\label{subsecApprox}
In order to obtain the approximation results we often impose a slight technical restriction (which can be removed, see for instance ~\cite{us}) on sequences of block partitions, which is defined as follows.
\begin{definition}[admissibility]
We say that a sequence $(\cR_N)_{N\geq 1}$ of block partitions of a block $R_0$ is \emph{admissible} if $\#(\cR_N) \leq N$ for all $N \geq 1$, and
\be
\label{defadmi}
\sup_{N\geq 1} \left( N^{\frac 1 d} \sup_{R \in \cR_N} \diam(R)\right) < \infty
\ee
\end{definition}

%

We recall that the approximation error is measured in $L_p$ norm, where the exponent $p$ is fixed and $1\leq p \leq \infty$.
We define $\tau\in (0, \infty)$ by
\be
\label{defTau}
\frac 1 \tau := \frac m d+ \frac 1 p.
\ee
In the following estimates we identified $d^m f(x)$ with an element of $\H_m$ according to
\be
\label{dmfHm}
\frac{d^m f(x)}{m!} \sim \sum_{|\balpha| = m} \frac{\partial^m f(x)}{\partial x^\balpha}  \frac{X^\balpha}{\balpha!}.
\ee
We now state the asymptotically sharp lower bound for the approximation error of a function $f$ on an admissible sequence of block partitions.
\begin{theorem}
\label{thLower}
Let $R_0$ be a block and let $f\in C^m(R_0)$. For any admissible sequence of block partitions $(\cR_N)_{N\geq 1}$ of $R_0$
$$
\liminf_{N\to \infty} N^{\frac m d}\|f-\interp_{\cR_N} f\|_{L_p(R_0)} \geq \left\|K_I\left(\frac{d^m f}{m!}\right)\right\|_{L_\tau(R_0)}.
$$
\end{theorem}

The next theorem provides an upper bound for the projection error of a function $f$ when an optimal sequence of block partitions is used. It confirms the sharpness of the previous theorem.

\begin{theorem}
\label{thUpper}
Let $R_0$ be a block and let $f\in C^m(R_0)$. Then there exists a (perhaps non-admissible) sequence $(\cR_N)_{N\geq 1}$, $\card\cR_N \leq N$, of block partitions of $R_0$ satisfying
\be
\label{upperEstim}
\limsup_{N\to \infty} N^{\frac m d}\|f-\interp_{\cR_N} f\|_{L_p(R_0)} \leq \left\|K_I\left(\frac{d^m f}{m!}\right)\right\|_{L_\tau(R_0)}.
\ee

Furthermore, for all $\ve>0$ there exists an admissible sequence $(\cR_N^\ve)_{N\geq 1}$ of block partitions of $R_0$ satisfying
\be
\label{upperEstimEps}
\limsup_{N\to \infty} N^{\frac m d}\|f-\interp_{\cR_N^\ve} f\|_{L_p(R_0)} \leq \left\|K_I\left(\frac{d^m f}{m!}\right)\right\|_{L_\tau(R_0)}+\ve.
\ee
\end{theorem}
An important feature of these estimates is the ``$\limsup$''.
Recall that the upper limit of a sequence $(u_N)_{N\geq N_0}$ is defined by
$$
\limsup_{N\to \infty} u_N := \lim_{N\to \infty} \sup_{n\geq N} u_n,
$$
and is in general strictly smaller than the supremum $\sup_{N\geq
N_0} u_N$. It is still an open question to find an appropriate upper
estimate of $\sup_{N\geq N_0} N^{\frac{m} d}  \|f-\interp_{\cR_N}
f\|_{L_p(R_0)}$ when optimally adapted block partitions are used.

In order to have more control of the quality of approximation on various parts of the domain we introduce a positive weight function $\Omega\in C^0(R_0)$.
For $1\leq p\leq \infty$ and for any $u\in L_p(R_0)$ as usual we define
\be 
\label{defWeight}
\|u\|_{L_p(R_0, \Omega)} := \|u\Omega\|_{L_p(R_0)}. 
\ee

\begin{remark}
\label{remarkWeight} Theorems \ref{thLower}, \ref{thUpper} and
\ref{thNoEps} below also hold when the norm $\| \cdot\|_{L_p(R_0)}$ (resp
$\| \cdot\|_{L_\tau(R_0)}$) is replaced with the weighted norm $\| \cdot\|_{L_p(R_0,
\Omega)}$ (resp $\| \cdot\|_{L_\tau(R_0, \Omega)}$) defined in \iref{defWeight}.
\end{remark}

In the following section we shall use some restrictive hypotheses on the interpolation operator in order to obtain an explicit formula for the shape function.
In particular, 
Propositions \ref{propOdd}, \ref{propEven}, and equation \iref{Kdmf} show that, under some assumptions, there exists a constant $C=C(\interp)>0$ such that
$$
\frac 1 {C} K_I\left(\frac{d^m f}{m!}\right) \leq \sqrt[d]{\left|\prod_{1\leq i\leq d} \frac{\partial^m f}{\partial x_i^m}\right|} \leq C K_I\left(\frac{d^m f}{m!}\right).
$$
These restrictive hypotheses also allow to improve slightly the estimate \iref{upperEstimEps} as follows.
\begin{theorem}
\label{thNoEps}
If the hypotheses of Proposition \ref{propOdd} or \ref{propEven} hold, and if
$
K_I\left(\frac{d^m f}{m!}\right) >0
$
everywhere on $R_0$, then there exists an \emph{admissible} sequence of partitions $(\cR_N)_{N\geq 1}$ which satisfies the optimal estimate \iref{upperEstim}.
\end{theorem}
The proofs of the Theorems \ref{thLower}, \ref{thUpper} and \ref{thNoEps} are given in \S\ref{proofTh}.
Each of these proofs can be adapted to weighted norms, hence establishing Remark \ref{remarkWeight}. Some details on how to adapt proofs for the case of weighted norms are provided at the end of each proof.

\section{Study of the error function}
\label{studyK}

In this section we perform  a close study of the error function $K_I$, since it plays a major role in our asymptotic error estimates. In the first subsection \S\ref{studyKGen} we investigate general properties which are valid for any continuous projection operator $\interp$.
However, we are not able to obtain an explicit form of $K_I$ under such general assumptions. 
Recall that in \S\ref{exampleI} we presented several possible choices of projection operators $\interp$ that seem more likely to be used in practice. In \S\ref{subsecproperties} we identify four important properties shared by these examples. 
These properties are used in \S \ref{subsecExact} to obtain an explicit form of $K_I$.

\subsection{General properties}
\label{studyKGen}
The error function $K$ obeys the following important invariance property with respect to diagonal changes of coordinates.
\begin{prop}
\label{propInv}  For all $\pi\in \H_m$ and all diagonal matrices $D$
with non-negative coefficients
$$
K_I(\pi\circ D) = (\det D)^{\frac m d} K_I(\pi).
$$
\end{prop}

\proof
We first assume that the diagonal matrix $D$ has positive diagonal coefficients. Let $\overline D$ be a diagonal matrix with positive diagonal coefficient and which satisfies $\det \overline D = 1$.
Let also $\pi\in \H_m$. Then
$$
\pi \circ (D \overline D) = \pi \circ ((\det D)^{\frac 1 d} \tilde D) = (\det D)^{\frac m d} \pi \circ \tilde D,
$$
where
$
\tilde D := (\det D)^{- \frac 1 d} D\overline D
$
satisfies $\det \tilde D = \det \overline D=1$ and is uniquely determined by $\overline D$.
According to \iref{KD} we therefore have
\begin{eqnarray*}
K_I(\pi \circ D) &=& \inf_{\substack{\det \overline D = 1\\ \overline D \geq 0}} \|\pi\circ (D \overline D) - \interp(\pi\circ (D \overline D))\|_{L_p(\mI^d)}\\
&=& (\det D)^{\frac m d} \inf_{\substack{\det \tilde D = 1\\ \tilde D \geq 0}} \|\pi\circ \tilde D - \interp(\pi\circ \tilde D)\|_{L_p(\mI^d)}\\
&=& (\det D)^{\frac m d} K_I(\pi),
\end{eqnarray*}
which concludes the proof in the case where $D$ has positive diagonal coefficients.\\
Let us now assume that $D$ is a diagonal matrix with non-negative diagonal coefficients and such that $\det(D) = 0$. 
Let $D'$ be a diagonal matrix with positive diagonal coefficients, and such that $D=DD'$ and $\det D' = 2$. We obtain 
$$
K_I(\pi \circ D) = K_I(\pi \circ (DD')) = 2^{\frac m d} K_I(\pi \circ D),
$$
which implies that $K_I(\pi \circ D) = 0$ and concludes the proof.
%
\sq

The next proposition shows that the exponent $p$ used for measuring the approximation error plays a rather minor role.
By $K_p$ we denote the error function associated with the exponent $p$.
\begin{prop}
\label{propEquiv} There exists a constant $c>0$ such that for all
$1\leq p_1\leq p_2\leq \infty$ we have on $\H_m$
$$
c K_\infty\leq K_{p_1} \leq K_{p_2} \leq K_{\infty}.
$$
\end{prop}

\proof 
For any function $f\in V = C^0(\mI^d)$ and for any $1\leq p_1\leq p_2\leq \infty$ by a standard convexity argument we obtain that
$$
\|f\|_{L_1(\mI^d)} \leq \|f\|_{L_{p_1}(\mI^d)} \leq \|f\|_{L_{p_2}(\mI^d)} \leq \|f\|_{L_{\infty}(\mI^d)}.
$$
Using \iref{KD}, it follows that
$$
K_1\leq K_{p_1}\leq K_{p_2}\leq K_\infty
$$
on $\H_m$.
Furthermore, the following semi norms on $\H_m$
$$
|\pi|_1 := \|\pi-\interp \pi\|_{L_1(\mI^d)} \text{ and } |\pi|_\infty := \|\pi-\interp \pi\|_{L_\infty(\mI^d)}
$$
vanish precisely on the same subspace of $\H_m$, namely $V_{\interp} \cap H_m = \{\pi \in
\H_m \sep \pi = \interp \pi\}$. Since $\H_m$ has finite dimension,
it follows that they are equivalent. Hence, there exists a constant
$c>0$ such that $c|\cdot |_\infty\leq |\cdot |_1$ on $\H_m$. Using
\iref{KD}, it follows that $cK_\infty\leq K_1$, which concludes the
proof. \sq

\subsection{Desirable properties of the projection operator}
\label{subsecproperties}
The examples of projection operators presented in \S\ref{exampleI} share some important properties which allow to obtain the explicit expression of the error function $K_I$. 
These properties are defined below and called $H_\pm$, $H_\sigma$, $H_*$ or $H_{**}$.
They are satisfied when operator $\interp$ is the interpolation at equispaced points (Definition \ref{interpEqui}), at Tchebychev points (Definition \ref{interpTche}), and usually on the most interesting sets of other points. They are also satisfied when $\interp$ is the $L_2(\mI^d)$ orthogonal projection onto $\P_k^*$ or $\P_k^{**}$ (Definition \ref{L2Proj}).

The first property reflects the fact that a coordinate $x_i$ on $\mI^d$ can be changed to $-x_i$, independently of the projection process.
\begin{definition}[$H_\pm$ hypothesis]
We say that the interpolation operator $\interp$ satisfies the $H_\pm$ hypothesis if for any diagonal matrix $D$ with entries in $\pm 1$ we have for all $f\in V$
$$
 \interp(f\circ D) = \interp(f) \circ D.
$$
\end{definition}

The next property implies that the different coordinates $x_1, \cdots, x_d$ on $\mI^d$ play symmetrical roles with respect to the projection operator.
\begin{definition}[$H_\sigma$ hypothesis]
If $M_\sigma$ is a permutation matrix, i.e. $(M_\sigma)_{ij} :=
\delta_{i \sigma(j)}$ for some permutation $\sigma$ of $\{1,\cdots,
d\}$, then for all $f\in V$
$$
\interp (f\circ M_\sigma) = \interp( f) \circ M_\sigma.
$$
\end{definition}

According to \iref{defk}, the projection operator $\interp$ reproduces the space of polynomials $\P_k$. However, in many situations the space $V_{\interp}$ of functions reproduced by $\interp$ is larger than $\P_k$. In particular $V_{\interp} = \P_k^{**}$ when $\interp$ is the interpolation on equispaced or Tchebychev points, and $V_{\interp} = \P_k$ (resp $\P_k^*$,  $\P_k^{**}$) when $\interp$ is the $L_2(\mI^d)$ orthogonal projection onto $\P_k$ (resp $\P_k^*$,  $\P_k^{**}$).

It is particularly useful to know whether the projection operator $\interp$ reproduces the elements of $\P_k^*$, and we therefore give a name to this property. Note that it clearly does not hold for the $L_2(\mI^d)$ orthogonal projection onto $\P_k$.
\begin{definition}[$H_*$ hypothesis]
The following inclusion holds :
$$
P_k^* \subset V_{\interp}.
$$
\end{definition}

On the contrary it is useful to know that some polynomials, and in particular pure powers $x_i^m$, are not reproduced by $\interp$.
\begin{definition}[$H_{**}$ hypothesis]
$$
\text{If } \ \sum_{1\leq i\leq d} \lambda_i x_i^m \in V_{\interp} \ \text{ then } \ (\lambda_1, \cdots, \lambda_d) = (0, \cdots , 0).
$$
\end{definition}
This condition obviously holds if $\interp(f)\in \P_k^{**}$ (polynomials of degree $\leq k$ in each variable) for all $f$. Hence, it holds for all the examples of projection operators given in the previous subsection \S\ref{exampleI}.

\subsection{Explicit formulas}
\label{subsecExact}

In this section
we provide the explicit expression for $K$ when some of the hypotheses $H_\pm$, $H_\sigma$, $H_*$ or $H_{**}$ hold.
Let $\pi\in \H_m$ and let $\lambda_i$ be the corresponding coefficient of $X_i^m$ in $\pi$, for all $1\leq i\leq d$.
We define
$$
K_*(\pi) := \sqrt[d]{\prod_{1\leq i\leq d} |\lambda_i|}
$$
and
$$
s(\pi) := \#\{1\leq i\leq d\sep \lambda_i >0\}.
$$
If $\frac{d^m f(x)}{m!}$ is identified by  \iref{dmfHm}  to an element of $\H_m$, then one has 
\be
\label{Kdmf}
K_*\left(\frac{d^m f(x)}{m!}\right) = \frac 1 {m!} \sqrt[d]{\left|\prod_{1\leq i\leq d} \frac{\partial^m f}{\partial x_i^m}(x)\right|}.
\ee

\begin{prop}
\label{propOdd}
If $m$ is odd and if $H_\pm$, $H_\sigma$ and $H_*$ hold, then
$$
K_p(\pi) = C(p) K_*(\pi),
$$
where
$$
 C(p) := \left\|\sum_{1\leq i\leq d} X_i^m - \interp \left(\sum_{1\leq i\leq d} X_i^m\right)\right\|_{L_p(\mI^d)}>0.
$$
\end{prop}

\begin{prop}
\label{propEven}
If $m$ is even and if $H_\sigma$, $H_*$ and $H_{**}$ hold then
$$
K_p(\pi) = C(p,s(\pi)) \, K_*(\pi).
$$
Furthermore,
\be
\label{eqEven}
C(p,0) = C(p,d) =
\left\|\sum_{1\leq i\leq d} X_i^m - \interp
\left(\sum_{1\leq i\leq d} X_i^m\right)\right\|_{L_p(\mI^d)}>0.
\ee
Other constants $C(p,s)$ are positive and obey $C(p,s) =
C(p,d-s)$.\\
\end{prop}

\noindent
Next we turn to the proofs of Propositions \ref{propOdd} and \ref{propEven}.
\paragraph{Proof of Proposition \ref{propOdd}}
Let $\pi\in \H_m$ and let $\lambda_i$ be the coefficient of $X_i^m$ in $\pi$.
Denote by
$$
\pi_* := \sum_{ 1\leq i\leq d} \lambda_i X_i^m
$$
so that $\pi-\pi_* \in \P_k^*$ and, more generally, $\pi \circ D - \pi_* \circ D \in \P_k^*$ for any diagonal matrix $D$.
The hypothesis $H_*$ states that the projection operator $\interp$ reproduces the elements of $\P_k^*$, and therefore 
$$
\pi\circ D - \interp (\pi\circ D) = \pi_*\circ D - \interp (\pi_*\circ D).
$$
Hence,
$
K_I(\pi) = K_I(\pi_*)
$
according to \iref{KD}.
If there exists $i_0$, $1\leq i_0\leq d$, such that $\lambda_{i_0} = 0$, then we denote by $D$ the diagonal matrix of entries $D_{ii} = 1$ if $i\neq i_0$ and $0$ if $i=i_0$. Applying Proposition \ref{propInv} we find 
$$
K_I(\pi) = K_I(\pi_*) = K_I(\pi_*\circ D) = (\det D)^{\frac m d} K_I(\pi_*) = 0.
$$
which concludes the proof.
We now assume that all the coefficients $\lambda_i$, $1\leq i\leq d$, are different from $0$, and we denote by $\ve_i$ be the sign of $\lambda_i$. Applying Proposition
\ref{propInv} to the diagonal matrix $D$ of entries $|\lambda_i|^{\frac
1 m}$ we find that
$$
K_I(\pi) = K_I(\pi_*) = (\det D)^{\frac m d} K_I(\pi_*\circ D^{-1}) = K_*(\pi) \, K_I\left(\sum_{1\leq i\leq d} \ve_i X_i^m\right).
$$
Using the $H_\pm$ hypothesis with the diagonal matrix $D$ of entries
$D_{ii} = \ve_i$, and recalling that $m$ is odd, we find that
$$
K_I\left(\sum_{1\leq i\leq d} \ve_i X_i^m\right) = K_I\left(\sum_{1\leq i\leq d}  X_i^m\right).
$$
We now define the functions
$$
g_i := X_i^m - \interp(X_i^m) \text{ for } 1\leq i \leq d.
$$
It follows from \iref{KD} that
$$
K_I\left(\sum_{1\leq i\leq d}  X_i^m\right) = \inf_{\prod_{1\leq
i\leq d} a_i = 1} \left\|\sum_{1\leq i\leq d} a_i
g_i\right\|_{L_p(\mI^d)},
$$
where the infimum is taken over all $d$-vectors of positive reals of product $1$.
Let us consider such a $d$-vector $(a_1, \cdots , a_d)$, and a permutation $\sigma$ of the set $\{1, \cdots,
d\}$.
The $H_\sigma$ hypothesis implies that the quantity
$$
\left\|\sum_{1\leq i\leq d} a_{\sigma(i)} g_i\right\|_{L_p(\mI^d)}
$$
is independent of $\sigma$. Hence, summing over all permutations, we obtain \be
\label{sumPerm} \left\|\sum_{1\leq i\leq d} a_i
g_i\right\|_{L_p(\mI^d)} = \frac 1 {d!} \sum_{\sigma}
\left\|\sum_{1\leq i\leq d} a_{\sigma(i)} g_i\right\|_{L_p(\mI^d)}
\geq \frac 1 d \left\|\sum_{1\leq i\leq d} g_i\right\|_{L_p(\mI^d)}
\sum_{1\leq i\leq d} a_i. \ee 
The right-hand side is minimal when $ a_1 =
\cdots = a_d = 1 $, which shows that
$$
\left\|\sum_{1\leq i\leq d} a_i g_i\right\|_{L_p(\mI^d)} \geq  \left\|\sum_{1\leq i\leq d} g_i\right\|_{L_p(\mI^d)} = C(p)
$$
with equality when $a_i=1$ for all $i$. 
Note as a corollary that 
\be
\label{existsRPos}
K_I(\pi_\ve) = \|\pi_\ve - \interp(\pi_\ve)\|_{L_p(\mI^d)} = C(p) \ \text{ where } \ \pi_\ve = \sum_{1\leq i\leq d} \ve_i X_i^m.
\ee
It remains to prove that $C(p)>0$. Using the hypothesis $H_\pm$, we find that for all $\mu_i \in\{\pm 1\}$ we have
$$
 \left\|\sum_{1\leq i\leq d}\mu_i g_i\right\|_{L_p(\mI^d)} = C(p).
$$
In particular, for any $1\leq i_0\leq d$ one has 
$$
2\|g_{i_0}\|_{L_p(\mI^d)} \leq  \left\|\sum_{1\leq i\leq d} g_i\right\|_{L_p(\mI^d)} +  \left\|2g_{i_0} - \sum_{1\leq i\leq d} g_i\right\|_{L_p(\mI^d)} \leq 2 C(p).
$$
If $C(p) = 0$, it follows
that $g_{i_0} = 0$ and therefore that $X_{i_0}^m = \interp(X_{i_0}^m)$, for any $1\leq
i_0\leq d$. Using the assumption $H_*$, we find that the projection operator
$\interp$ reproduces all the polynomials of degree $m= k+1$, which
contradicts the definition \iref{defk} of the integer $k$.

\sq

\paragraph{Proof of proposition \ref{propEven}}
We define $\lambda_i$, $\pi_*$ and $\ve_i\in\{\pm 1\}$ as before and
we find, using similar reasoning, that
$$
K_I(\pi) = K_*(\pi) K_I\left(\sum_{1\leq i\leq d} \ve_i X_i^m\right).
$$
\noindent
For $1\leq s\leq d$ we define
$$
C(p, s) := K_I\left(\sum_{1\leq i\leq s}  X_i^m - \sum_{s+1\leq i\leq d}  X_i^m\right).
$$
From the hypothesis $H_\sigma$ it follows that $K_I(\pi) = K_*(\pi)
C(p,s(\pi))$.

Using again $H_\sigma$ and the fact that $K_I(\pi) = K_I(-\pi)$ for all $\pi \in \H_m$, we find that
$$
C(p,s) = K_I\left(\sum_{1\leq i\leq s}  X_i^m - \sum_{s+1\leq i\leq d}  X_i^m\right) = K_I\left(-\left(\sum_{1\leq i\leq d-s}  X_i^m - \sum_{d-s+1\leq i\leq d}  X_i^m\right)\right) = C(p,d-s).
$$

We define $g_i := X_i^m - \interp(X_i^m)$, as in the proof of Proposition \ref{propOdd}.
We obtain the expression for $C(p,0)$ by summing over all permutations as in  \iref{sumPerm}
$$
C(p,0) = \left\|\sum_{1\leq i\leq d} g_i\right\|_{L_p(\mI^d)}.
$$
This concludes the proof of the first part of Proposition \ref{propEven}.
We now prove that $C(p,s)>0$ for all $1\leq p\leq \infty$ and all $s\in \{0, \cdots, d\}$.
To this end we define the following quantity on $\R^d$
$$
 \|a\|_K := \left\|\sum_{1\leq i\leq d}a_i  g_i\right\|_{L_p(\mI^d)}.
$$
Note that $\|a\|_K = 0$ if and only if
$$
\sum_{1\leq i\leq d} a_i X_i^m = \sum_{1\leq i\leq d} a_i \interp(X_i^m),
$$
and the hypothesis $H_{**}$ precisely states that this equality occurs if and only if $a_i = 0$, for all $1\leq i\leq d$. Hence, $\|\cdot \|_K$ is a norm on $\R^d$.
Furthermore, let
$$
E_s := \left\{a\in \R_+^s\times \R_-^{d-s}\sep \prod_{1\leq i \leq
d} |a_i| = 1\right\}
$$
Then
$$
C(p,s) = \inf_{a\in E_s} \|a\|_K.
$$
Since $E_s$ is a closed subset of $\R^d$, which does not contain the
origin, this infimum is attained. It follows that $C(p,s)>0$, and that there exists a rectangle $R_\ve$ of unit volume such that
\be
\label{existsREven}
K_I(\pi_\ve) = \|\pi_\ve - \interp \pi_\ve\|_{L_p(R_\ve)} = C(p, s(\pi_\ve)) \ \text{ where } \ \pi_\ve = \sum_{1\leq i\leq d} \ve_i X_i^m.
\ee
\sq

\section{Proof of the approximation results}
\label{proofTh}

In this section, let the
block $R_0$, the integer $m$, the function $f\in C^m(R_0)$ and the exponent $p$ be fixed. We conduct our proofs for $1\leq p<\infty$ and provide comments on how to adjust our arguments for the case $p=\infty$.

For each $x\in \R_0$ by $\mu_x$ we denote
 the  $m$-th degree Taylor
polynomial of $f$ at the point $x$
\be
\label{defmux}
\mu_x = \mu_x(X) := \sum_{|\balpha| \leq m}\frac{\partial^m
f}{\partial x^\balpha}(x) \, \frac{(X-x)^\balpha}{\balpha!},
\ee
and we define $\pi_x$ to be the homogeneous component of degree $m$ in
$\mu_x$,
\be
\label{defpix}
\pi_x = \pi_x(X) :=  \sum_{|\balpha| = m}\frac{\partial^m f}{\partial x^\balpha}(x) \, \frac{X^\balpha}{\balpha!}.
\ee
Since $\pi_x$ and $\mu_x$ are polynomials of degree $m$, their $m$-th derivative is constant, and clearly $d^m \pi_x=d^m \mu_x = d^m
f(x)$. 
In particular, for any $x\in R_0$ the polynomial $\mu_x - \pi_x$ belongs to $\P_k$ (recall that $k=m-1$) and is therefore reproduced by the projection operator $\interp$. It follows that for
any $x\in R_0$ and any block $R$
\be
\label{muzPiz} \pi_x- \interp_R \pi_x = \mu_x
-\interp_R \mu_x.
\ee
In addition, we introduce a measure $\rho$ of
the degeneracy of a block $R$
$$
\rho(R) := \frac{\diam(R)^d}{|R|}.
$$
Given any function $g\in C^m(R)$ and any $x\in R$ we can define, similarly to \iref{defpix}, a polynomial $\hat \pi_x\in \H_m$ associated to $g$ at $x$.
We then define 
\be
\label{normdmg}
\|d^m g\|_{L_\infty(R)}:=\sup_{x\in R} \left(\sup_{|u| = 1} |\hat \pi_x(u)|\right).
\ee
\begin{prop}
There exists a constant $C = C(m,d)>0$ such that for any block $R$ and any function $g\in C^m(R)$
\be
\label{localIso}
\|g-\interp_R g\|_{L_p(R)} \leq C |R|^{\frac 1 \tau} \rho(R)^{\frac m d} \|d^m g\|_{L_\infty(R)}.
\ee
\end{prop}

\proof
Let $x_0\in R$ and let $g_0$ be the Taylor polynomial for $g$ of degree $m-1$ at point $x_0$ which is defined as follows
$$
g_0(X) := \sum_{|\balpha| \leq m-1}\frac{\partial^\balpha f(x_0)}{\partial x^\balpha}\frac{(X-x_0)^\balpha}{\balpha!}.
$$
Let $x\in R$ and let $x(t) =  x_0+ t (x-x_0)$.
We have
$$
g(x) = g_0(x)+ \int_{t=0}^1 \frac{d^m g_{x(t)}(x-x_0)}{m!} (1-t)^m dt.
$$
Hence,
\be
\label{ineqGG0}
|g(x) - g_0(x)|\leq  \int_{t=0}^1 \|d^m g\|_{L_\infty(R)} |x-x_0|^m (1-t)^mdt \leq \frac 1 {m+1} \|d^m g\|_{L_\infty(R)} \diam(R)^m.
\ee
Since $g_0$ is a polynomial of degree at most $m-1$, we have $g_0 = \interp g_0$. Hence,
\begin{eqnarray*}
\|g-\interp_R g\|_{L_p(R)} &\leq &  |R|^{\frac 1 p} \|g-\interp_R g\|_{L_\infty(R)}\\
& =& |R|^{\frac 1 p} \|(g-g_0)-\interp_R (g-g_0)\|_{L_\infty(R)}\\
&\leq & (1+C_{\interp}) |R|^{\frac 1 p} \|g-g_0\|_{L_\infty(R)},
\end{eqnarray*}
where $C_{\interp}$ is the operator norm of $\interp : V\to V$. Combining this estimate with \iref{ineqGG0}, we obtain \iref{localIso}.

\sq

\subsection{Proof of Theorem \ref{thLower} (Lower bound)}

The following lemma allows us to bound the interpolation error of $f$ on the block $R$ from below.
\begin{lemma}
\label{lemmaLower} For any block $R\subset R_0$ and $x\in R$ we have
$$
\|f- \interp_R f\|_{L_p(R)} \geq |R|^{\frac 1 \tau} \left(K_I(\pi_x) - \omega(\diam R) \rho(R)^{\frac m d}\right),
$$
where the function $\omega$ is positive, depends only on $f$ and $m$, and satisfies $\omega(\delta) \to 0$ as $\delta\to 0$.
\end{lemma}

\proof
Let $h:= f-\mu_x$, where $\mu_x$ is defined in \iref{defmux} Using \iref{muzPiz}, we obtain
\begin{eqnarray*}
\|f- \interp_R f\|_{L_p(R)} & \geq & \|\pi_x - \interp_R \pi_x\|_{L_p(R)} - \|h - \interp_R h\|_{L_p(R)}\\
& \geq &  |R|^{\frac 1 \tau} K_I(\pi_x) - \|h - \interp_R
h\|_{L_p(R)},
\end{eqnarray*}
and according to \iref{localIso} we have
$$
 \|h - \interp_R h\|_{L_p(R)} 
\leq C_0 |R|^{\frac 1 \tau}
\rho(R)^{\frac m d} \|d^m h\|_{L_\infty(R)}.
$$
Observe that
$$
\|d^m h\|_{L_\infty(R)} = \|d^m f - d^m\pi_x\|_{L_\infty(R)}  =
\|d^m f - d^m f(x)\|_{L_\infty(R)}. 
$$
We introduce the modulus of continuity $\omega_*$ of the $m$-th derivatives of $f$. 
\be
\label{defOmega}
\omega_*(r) := \sup_{\substack{x_1,x_2\in R_0 :\\ |x_1 - x_2|\leq r}} \|d^m f(x_1)-d^m f(x_2)\|= \sup_{\substack{x_1,x_2\in R_0 :\\ |x_1-x_2|\leq r}} \left( \sup_{|u|\leq 1} |\pi_{x_1}(u) - \pi_{x_2} (u)|\right)
\ee
By setting $\omega = C_0\, \omega_*$
we conclude the proof of
this lemma.
\sq

We now consider an admissible sequence of block partitions $(\cR_N)_{N\geq 0}$. For all $N\geq 0$,
$R\in \cR_N$ and $x\in R$, we define
$$
\phi_N(x) := |R| \qquad \hbox{and} \qquad
\psi_N(x) := \left(K_I(\pi_x) - \omega(\diam(R)) \rho(R)^{\frac m d} \right)_+,
$$
where $\lambda_+ := \max\{\lambda,0\}$.
We now apply Holder's inequality
$
\int_{R_0} f_1 f_2 \leq \|f_1\|_{L_{p_1}(R_0)}\|f_2\|_{L_{p_2}(R_0)}
$
with the functions
$$
f_1 = \phi_N^{\frac {m \tau} d} \psi_N^\tau \ \text{ and } f_2 = \phi_N^{-\frac {m \tau} d}
$$
and the exponents
$
p_1 = \frac p \tau \ \text{ and } \ p_2 = \frac d {m\tau}.
$
Note that $\frac 1 {p_1}+ \frac 1 {p_2} =  \tau \left(\frac 1 p + \frac m d\right) = 1$. Hence,
\be
\label{holderPsi}
\int_{R_0} \psi_N^\tau \leq \left(\int_{R_0} \phi_N^{\frac{m p} d}
\psi_N^{p}\right)^{\frac \tau p} \left(\int_{R_0}
\phi_N^{-1}\right)^{\frac {m\tau} d}.
\ee

\noindent
Note that $\int_{R_0} \phi_N^{-1} = \# (\cR_N)\leq N$. Furthermore, if
$R\in \cR_N$ and $x\in R$ then according to Lemma \ref{lemmaLower}
$$
 \phi_N(x)^{\frac m d} \psi_N(x)
=|R|^{\frac 1 \tau-\frac 1 p}  \psi_N(x) \leq|R|^{-\frac 1 p} \|f-\interp_R f\|_{L_p(R)}.
$$
Hence,
\be
\label{intphipsi}
\left[\int_{R_0} \phi_N^{\frac{m p} d} \psi_N^{p}\right]^{\frac 1 p} \leq \left[\sum_{R\in \cR_N} \frac 1 {|R|}  \int_R  \|f-\interp_R f\|_{L_p(R)}^p\right]^{\frac 1 p}=
\|f-\interp_R f\|_{L_p(R_0)}.
\ee
Inequality \iref{holderPsi} therefore leads to
\be
\label{upperPsi}
 \|\psi_N\|_{L_\tau(R_0)} \leq \|f- \interp_{\cR_N} f\|_{L_p(R_0)} N^{\frac m d}.
\ee

\noindent
Since the sequence $\seqR$ is admissible, there exists a constant $C_A>0$ such that for all $N$ and all $R\in \cR_N$ we have $\diam(R)\leq C_AN^{-\frac 1 d}$.
We introduce a subset of $\cR'_N\subset \cR_N$ which collects the most degenerate blocks
$$
\cR'_N = \{ R\in \cR_N \sep \rho(R)\geq \omega(C_AN^{-\frac 1
d})^{-\frac{1} m}\},
$$
where $\omega$ is the function defined in Lemma \ref{lemmaLower}. 
By $R'_N$ we denote the portion of $R_0$ covered by $\cR'_N$.
For all $x\in R_0\bs R'_N$
we obtain
$$
\psi_N(x)\geq K_I(\pi_x) -\omega(C_A N^{-\frac 1 d})^{1-\frac 1 d}.
$$
We define $\ve_N := \omega(C_A N^{-\frac 1 d})^{1-\frac 1 d}$ and we notice that $\ve_N \to 0$ as $N \to \infty$.
Hence,
$$
\begin{array}{ll}
 \|\psi_N\|_{L_\tau(R_0)}^\tau & \geq \left \|\left(K_I(\pi_x) -\ve_N\right)_+\right\|_{L_\tau(R_0\bs R'_N)}^\tau\\
 & \geq
 \left \|\left(K_I(\pi_x) -\ve_N\right)_+\right\|_{L_\tau(R_0)}^\tau
 -C^\tau |R'_N|,
 \end{array}
 $$
where $C:=\max_{x\in R_0}K_I(\pi_x)$. Next we observe
that $|R'_N|\to 0$ as $N\to +\infty$: indeed
for all $R\in \cR'_N$ we have
$$
|R| = \diam(R)^d \rho(R)^{-1} \leq C_A^d N^{-1} \omega(C_A N^{-\frac 1 d})^{\frac 1 m}.
$$
Since $\card(\cR'_N)\leq N$, we obtain $|R'_N|\leq C_A^d \omega(C_A N^{-\frac 1 d})^{\frac 1 m}$, and the right-hand side tends to $0$ as $N\to \infty$. We thus obtain
$$
\liminf_{N\to \infty} \|\psi_N\|_{L_\tau(R_0)}
\geq \lim_{N\to \infty}   \left \|\left(K_I(\pi_x) -\ve_N\right)_+\right\|_{L_\tau(R_0)}
 = \|K_I(\pi_x)\|_{L_\tau(R_0)}.
$$
Combining this result with \iref{upperPsi}, we conclude the proof of the announced estimate.

Note that this proof also works with the exponent $p = \infty$ by changing
$$
\left(\int_{R_0} \phi_N^{\frac{m p} d}
\psi_N^{p}\right)^{\frac \tau p} \ \text{ into } \ \|\phi_N^{\frac m d}
\psi_N\|_{L_\infty(R_0)}^\tau
$$
in \iref{holderPsi} and performing the standard modification in \iref{intphipsi}.
\begin{remark}
As announced in Remark \ref{remarkWeight}, this proof can be adapted to the weighted norm $\|\cdot\|_{L_p(R_0, \Omega)}$ associated to a positive weight function $\Omega\in C^0(R_0)$ and defined in \iref{defWeight}. For that purpose let $r_N := \sup \{ \diam(R) \sep R \in \cR_N\}$ and let
$$
\Omega_N(x) := \inf_{\substack{x'\in R_0\\ |x-x'|\leq r_N}} \Omega(x').
$$
The sequence of functions $\Omega_N$ increases with $N$ and tends uniformly to $\Omega$ as $N\to \infty$.
If $R\in \cR_N$ and $x\in R$, then
$$
\|f-\interp_R f\|_{L_p(R,\Omega)}
\geq  \Omega_N(x) \|f-\interp_R f\|_{L_p(R)}.
$$
The main change in the proof is that the function $\psi_N$ should be replaced with $\psi'_N := \Omega_N \psi_N$. Other details are left to the reader.
\end{remark}
\sq

\subsection{Proof of the upper estimates}

The proof of Theorems \ref{thUpper} (and \ref{thNoEps}) is based on the actual construction of an asymptotically optimal sequence of block partitions. To that end we introduce the notion of a local block specification.

\begin{definition}{\bf{(local block specification)}}
\label{defBlockSpec}
A local block specification on a block $R_0$ is a (possibly discontinuous) map $x \mapsto R(x)$ which associates to each point $x\in R_0$ a block $R(x)$, 
and such that
\begin{itemize}
\item
The volume
$
|R(x)|
$
is a positive continuous function of the variable $x \in R_0$.
\item The diameter is bounded : $\sup \{\diam(R(x))\sep x \in R_0\}<\infty$. 
\end{itemize}
\end{definition}

The following lemma shows that it is possible to build sequences of block partitions of $R_0$ adapted in a certain sense to a local block specification.

\begin{lemma}
\label{lemmaSeqBlock}
Let $R_0$ be a block in $\RR^d$ and let $x\mapsto R(x)$ be a local block specification on $R_0$. Then there exists a sequence $(\cP_n)_{n\geq 1}$ of block partitions of $R_0$,
$
\cP_n = \cP_n^1 \cup \cP_n^2,
$
satisfying the following properties.
\begin{itemize}
\item (The number of blocks in $\cP_n$ is asymptotically controlled)
\be \label{limCardRn} \lim_{n \to \infty} \frac{\#(\cP_n)}{n^{2d}} =
\int_{R_0} |R(x)|^{-1} dx. \ee
\item (The elements of $\cP_n^1$ follow the block specifications)
For each $R\in \cP_n^1$ there exists $y\in R_0$ such that
\be
\label{n2Ry}
R \text{ is a translate of } n^{-2} R(y), \text{ and } |x-y| \leq \frac{\diam(R_0)} n \text{ for all } x\in R.
\ee
\item (The elements of $\cP_n^2$ have a small diameter)
\be \label{smallDiam} \lim_{n \to \infty } \left( n^2 \sup_{R\in
\cP_n^2} \diam(R)\right) =0. \ee
\end{itemize}
\end{lemma}

\proof
See Appendix.
\sq

We recall that the block $R_0$, the exponent $p$ and the function $f\in C^m(R_0)$ are fixed, and that at each point $x\in R_0$ the polynomial $\pi_x\in \H_m$ is defined by \iref{defpix}.
The sequence of block partitions described in the previous lemma is now used to obtain an asymptotical error estimate. 
\begin{lemma}
\label{lemmaNn} Let $x \mapsto R(x)$ be a local block specification
such that for all $x\in R_0$ \be \label{unitError} \|\pi_x
-\interp_{R(x)}(\pi_x)\|_{L_p(R(x))} \leq 1. \ee Let $(\cP_n)_{n\geq
1}$ be a sequence of block partitions satisfying the properties of
Lemma \ref{lemmaSeqBlock}, and let for all $N\geq 0$
$$
n(N) := \max\{ n\geq 1  \sep \# (\cP_n) \leq N\}.
$$
Then $\cR_N := \cP_{n(N)}$ is an admissible sequence of block partitions and
\be
\label{limsupR}
\limsup_{N\to \infty} N^{\frac m d} \|f-\interp_{\cR_N} f \|_{L_p(R_0)} \leq \left(\int_{R_0} R(x)^{-1} dx\right)^{\frac 1 \tau}.
\ee
\end{lemma}

\proof
Let $n \geq 0$ and let $R\in \cP_n$.
If $R\in \cP_n^1$ then let $y\in R_0$ be as in \iref{n2Ry}. Using \iref{localIso} we find
\begin{eqnarray*}
\|f-\interp_R f\|_{L_p(R)} &\leq& \|\pi_y-\interp_R \pi_y\|_{L_p(R)}  + \|(f-\pi_y)-\interp_R (f-\pi_y)\|_{L_p(R)} \\
&\leq & n^{-\frac {2d} \tau}  \|\pi_y-\interp_{R(y)} \pi_y\|_{L_p(R(y))} + C |R|^{\frac 1 p} \diam(R)^m \|d^m f-d^m \pi_y\|_{L_{\infty(R)}}\\
&\leq &  n^{-\frac {2d} \tau}  + C n^{-\frac {2d} \tau} |R(y)|^{\frac 1 p} \diam(R(y))^m \|d^m f-d^m f(y)\|_{L_\infty(R)}\\
&\leq & n^{-\frac {2d} \tau}  (1+C' \omega_*(n^{-1}\diam(R_0))),
\end{eqnarray*}
where we defined $C' := C \sup_{y\in \R_0} |R(y)|^{\frac 1 p} \diam(R(y))^m$, which is finite by Definition \ref{defBlockSpec}. 
We denoted by $\omega_*$ the modulus of continuity of the $m$-th derivatives of $f$ which is defined at \iref{defOmega}.
We now define for all $n \geq 1$, 
$$
\delta_n := n^2 \sup_{R \in \cP_n^2} \diam(R).
$$
According to \iref{smallDiam} one has $\delta_n \to 0$ as $n \to \infty$.
If $R\in \cP_n^2$, then $\diam(R)\leq n^{-2} \delta_n$ and therefore $|R|\leq \diam(R)^d\leq n^{-2d} \delta_n^d$. Using again \iref{localIso}, and recalling that $\frac 1 \tau = \frac m d+ \frac 1 p$ we find
$$
\|f-\interp_R f\|_{L_p(R)} \leq C |R|^{\frac 1 p} \diam(R)^m \|d^m f\|_{L_{\infty(R_0)}} \leq C'' n^{-\frac {2d} \tau} \delta_n^{\frac d \tau}
$$
where $C'' = C \|d^m f\|_{L_{\infty(R_0)}}$.
From the previous observations it follows that
$$
\|f-\interp_{\cP_n} f\|_{L_p(R_0)} \leq \#(\cP_n)^{\frac 1 p}
\max_{R\in \cP_n} \|f-\interp_R f\|_{L_p(R)} \leq  \#(\cP_n)^{\frac
1 p}  n^{-\frac {2d} \tau}  \max\{1+
C'\omega_*(n^{-1}\diam(R_0)), \, C'' \delta_n^{\frac d
\tau}\}.
$$
Hence,
$$
\limsup_{n\to \infty} \#(\cP_n)^{-\frac 1 p}  n^{\frac {2d} \tau}\|f-\interp_{\cP_n} f\|_{L_p(R_0)}  \leq 1.
$$
Combining the last equation with \iref{limCardRn}, we obtain
$$
\limsup_{n\to \infty} \#(\cP_n)^{\frac m d} \|f-\interp_{\cP_n} f\|_{L_p(R_0)}  \leq \left(\int_{R_0} R(x)^{-1} dx\right)^{\frac 1 \tau}.
$$
The sequence of block partitions $\cR_N := \cP_{n(N)}$
clearly satisfies $\#(\cR_N)/N\to 1$ as $N
\to \infty$ and therefore leads to the announced equation \iref{limsupR}.
Furthermore, it follows from the boundedness of $\diam(R(x))$ on
$R_0$ and the properties of $\cP_n$ described in Lemma
\ref{lemmaSeqBlock} that
$$
\sup_{n\geq 1}\left( \#(\cP_n)^{\frac 1 d} \sup_{R \in \cP_n}
\diam(R) \right)< \infty
$$
which implies that $\cR_N$ is an admissible sequence of partitions.
\sq

We now choose adequate local block specifications in order to obtain
the estimates announced in Theorems \ref{thUpper} and \ref{thNoEps}.
For any $M\geq \diam(\mI^d) = \sqrt d$ 
we define the modified error function \be \label{defKM}
K_M(\pi) := \inf_{\substack{|R| = 1,\\ \diam(R)\leq M}} \|\pi - \interp_R \pi\|_{L_p(R)}, 
\ee
where the infimum is taken on blocks of unit volume and diameter smaller that $M$.
It follows from a compactness argument that this infimum is attained and that $K_M$ is a continuous function on $\H_m$. Furthermore, for all $\pi\in \H_m$, $M\mapsto K_M(\pi)$ is a decreasing function of $M$ which tends to $K_I(\pi)$ as $M \to \infty$.

For all $x\in R_0$ we denote by $R_M^*(x)$ a block which realises the infimum in $K_M(\pi_x)$. Hence,
$$
|R_M^*(x)| = 1, \ \diam(R_M^*(x))\leq M, \text{ and } K_M(\pi_x) = \|\pi_x - \interp_{R_M^*(x)} \pi_x\|_{L_p(R_M^*(x))}
$$
We define a local block specification on $R_0$ as follows
\be
\label{defRM}
R_M(x) := (K_M(\pi_x) + M^{-1})^{- \frac \tau d} R_M^*(x).
\ee
We now observe that
$$
 \|\pi_x - \interp_{R_M(x)} \pi_x\|_{L_p(R_M(x))}  = K_M(\pi_x) (K_M(\pi_x)+ M^{-1})^{-1} \leq 1.
$$
Hence, according to Lemma \ref{lemmaNn}, there exists a sequence $(\cR_N^M)_{N\geq 1}$ of block partitions of $R_0$ such that
$$
\limsup_{N\to \infty} N^{\frac m d} \|f-\interp_{\cR^M_N} f \|_{L_p(R_0)} \leq \|K_M(\pi_x)+M^{-1}\|_{L_\tau (R_0)}.
$$
Using our previous observations on the function $K_M$, we see that
$$
\lim_{M \to \infty} \|K_M(\pi_x)+M^{-1}\|_{L_\tau (R_0)} = \|K_I(\pi_x)\|_{L_\tau (R_0)}.
$$
Hence, given $\ve >0$ we can choose $M(\ve)$ large enough in such a way that
$$
\|K_{M(\ve)}(\pi_x)+M(\ve)^{-1}\|_{L_\tau (R_0)} \leq \|K_I(\pi_x)\|_{L_\tau (R_0)}+ \ve,
$$
which concludes the proof of the estimate \iref{upperEstimEps} of Theorem \ref{thUpper}.

For each $N$ let $M=M(N)$ be such that
$$
N^{\frac m d} \|f-\interp_{\cR^M_N} f \|_{L_p(R_0)} \leq \|K_M(\pi_x)+M^{-1}\|_{L_\tau (R_0)}+M^{-1}
$$
and $M(N) \to \infty$ as $N \to \infty$. Then the (perhaps non admissible) sequence of block partitions $\cR_N := \cR_N^{M(N)}$ satisfies \iref{upperEstim} which concludes the proof of Theorem \ref{thUpper}.
\sq

We now turn to the proof of Theorem \ref{thNoEps}, which follows the same scheme for the most. 
There exists $d$ functions $\lambda_1(x)
, \cdots, \lambda_d(x) \in C^0(R_0)$, and a function $x \mapsto \pi_*(x) \in \P^*_k$ such that
for all $x\in R_0$ we have
$$
\pi_x = \sum_{1\leq i\leq d} \lambda_i(x) X_i^m + \pi_*(x).
$$
The hypotheses of Theorem \ref{thNoEps} state that $K_I\left(\frac {d^m f}{m!}\right) = K_I(\pi_x)$ does not vanish on $R_0$. It follows from Propositions \ref{propOdd} and \ref{propEven} that the product $\lambda_1(x)\cdots \lambda_d(x)$ is nonzero for all $x\in R_0$.
We denote by $\ve_i\in \{\pm 1\}$ the sign of $\lambda_i$, which is therefore constant over the block $R_0$, and we define
$$
\pi_\ve := \sum_{1\leq i\leq d} \ve_i X_i^m
$$
The proofs of Propositions \ref{propEven} and \ref{propOdd} show that there exists a block $R_\ve$, satisfying $|R_\ve| = 1$, and such that $K_I(\pi_\ve) = \|\pi-\interp_{R_\ve} \pi\|_{L_p(R_\ve)}$.
By $D(x)$ we denote the diagonal matrix of entries $|\lambda_1(x)|, \cdots , |\lambda_d(x)|$, and we define
$$
R^*(x) := (\det D(x))^{\frac 1 {md} } D(x)^{- \frac 1 m} R_\ve.
$$
Clearly, $|R^*(x)| = 1$. Using \iref{changeRect} and the homogeneity
of $\pi_x\in \H_m$, we find that
$$
 \|\pi_x - \interp_{R^*(x)} \pi_x\|_{L_p(R^*(x))} = (\det D(x))^{\frac 1 d} K_I(\pi_\ve) = K_I(\pi_x).
$$
We then define the local block specification
\be
\label{defR}
R(x) := K_I(\pi_x)^{-\frac \tau d}R^*(x).
\ee
The admissible sequence $(\cR_N)_{N \geq 1}$ of block partitions constructed in Lemma \ref{lemmaNn} then satisfies the optimal upper estimate \iref{upperEstim}, which concludes the proof of Theorem \ref{thNoEps}. \sq

\begin{remark}[Adaptation to weighted norms]
Lemma \ref{lemmaNn} also holds if \iref{unitError} is replaced with
$$\Omega(x) \|\pi_x -\interp_{R(x)}(\pi_x)\|_{L_p(R(x))} \leq 1$$
and if the $L_p(R_0)$
norm is replaced with the weighted $L_p(R_0, \Omega)$ norm in \iref{limsupR}.
Replacing the block $R_M(x)$ defined in \iref{defRM} with
$$
R'_M(x) := \Omega(x)^{- \frac \tau d}R_M(x),
$$
one can easily obtain the extension of Theorem \ref{thUpper} to
weighted norms. Similarly, replacing $R(x)$ defined in \iref{defR}
with $R'(x) := \Omega(x)^{- \frac \tau d}R(x)$, one obtains the
extension of Theorem \ref{thNoEps} to weighted norms.
\end{remark}

\appendix
\begin{center}
\large APPENDIX
\end{center}

\section{Proof of Lemma \ref{lemmaSeqBlock}}
By $\cQ_n$ we denote  the standard partition of $R_0\in \RR^d$ in $n^d$ identical blocks of diameter $n^{-1} \diam(R_0)$ illustrated on the left in Figure 1.
For each $Q \in \cQ_n$ by $x_Q$ we denote the barycenter of $Q$ and we consider the tiling $\cT_Q$ of $\R^d$ formed with the block $n^{-2}R(x_Q)$ and its translates.
We define $\cP_n^1(Q)$ and $\cP_n^1$ as follows
$$
\cP_n^1(Q) := \{R\in \cT_Q \sep R \subset Q\} \ \text{ and } \ \cP_n^1 := \bigcup_{Q\in \cQ_n} \cP_n^1(Q).
$$
Comparing the areas, we obtain
$$
\card (\cP_n^1) = \sum_{Q\in \cQ_n} \cP_n^1(Q) \leq \sum_{Q \in \cQ_n} \frac{|Q|}{|n^{-2} R(x_Q)|} = n^{2 d} \sum_{Q \in \cQ_n} |Q| |R(x_Q)|^{-1}.
$$
From this point, using the continuity of $x \mapsto |R(x)|$, one can easily show that
$
\frac {\# (\cP_n^1)}{n^{2d}} \to \int_{R_0} |R(x)|^{-1} dx
$
as $n \to \infty$. Furthermore, the property \iref{n2Ry} clearly holds. In order to construct $\cP_n^2$, we first define two sets of blocks $\cP_n^{2*}(Q)$ and $\cP_n^{2*}$ as follows
$$
\cP_n^{2*}(Q) := \{R\cap Q \sep R \in \cT_Q \text{ and }R\cap \partial Q\neq \emptyset\} \ \text{ and } \ \cP_n^{2*} := \bigcup_{Q\in \cQ_n} \cP_n^{2*}(Q).
$$
Comparing the surface of $\partial Q$ with the dimensions of $R(x_Q)$, we find that
$$
\#(\cP_n^{2*}(Q)) \leq C n^{d-1}
$$
where $C$ is independent of $n$ and of $Q\in \cQ_n$.
Therefore, $\#(\cP_n^{2*})\leq C n^{2d-1}$. The set of blocks $\cP_n^2$ is then obtained by subdividing each block of $\cP_n^{2*}$ into $o(n)$ (for instance, $\lfloor \ln(n)\rfloor^d$) identical sub-blocks, in such a way that $\#(\cP_n^2)$ is $o(n^{2d})$ and that the requirement \iref{smallDiam} is met.

\begin{figure}
    \centering
    \includegraphics[width=4cm,height=4cm]{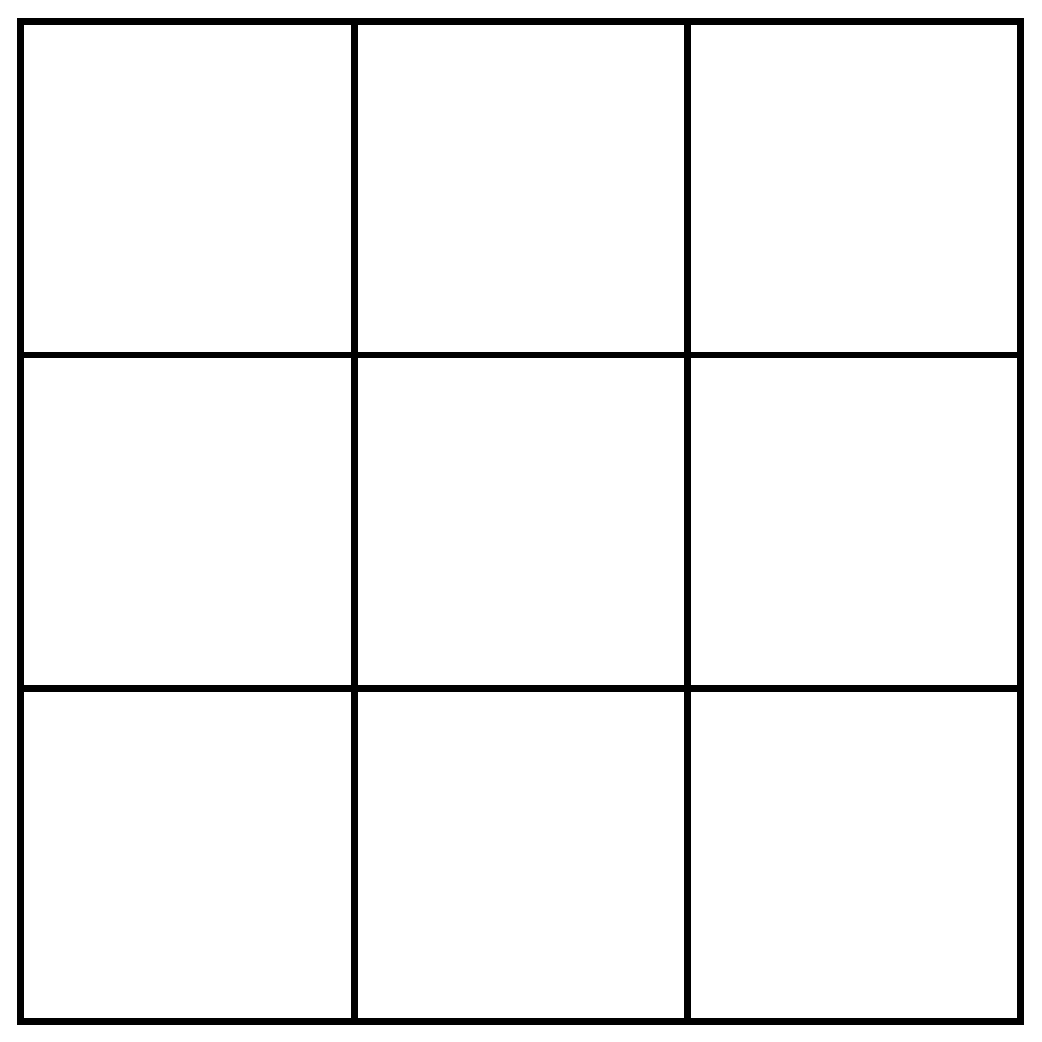}
    \hspace{1cm}
    \includegraphics[width=4cm,height=4cm]{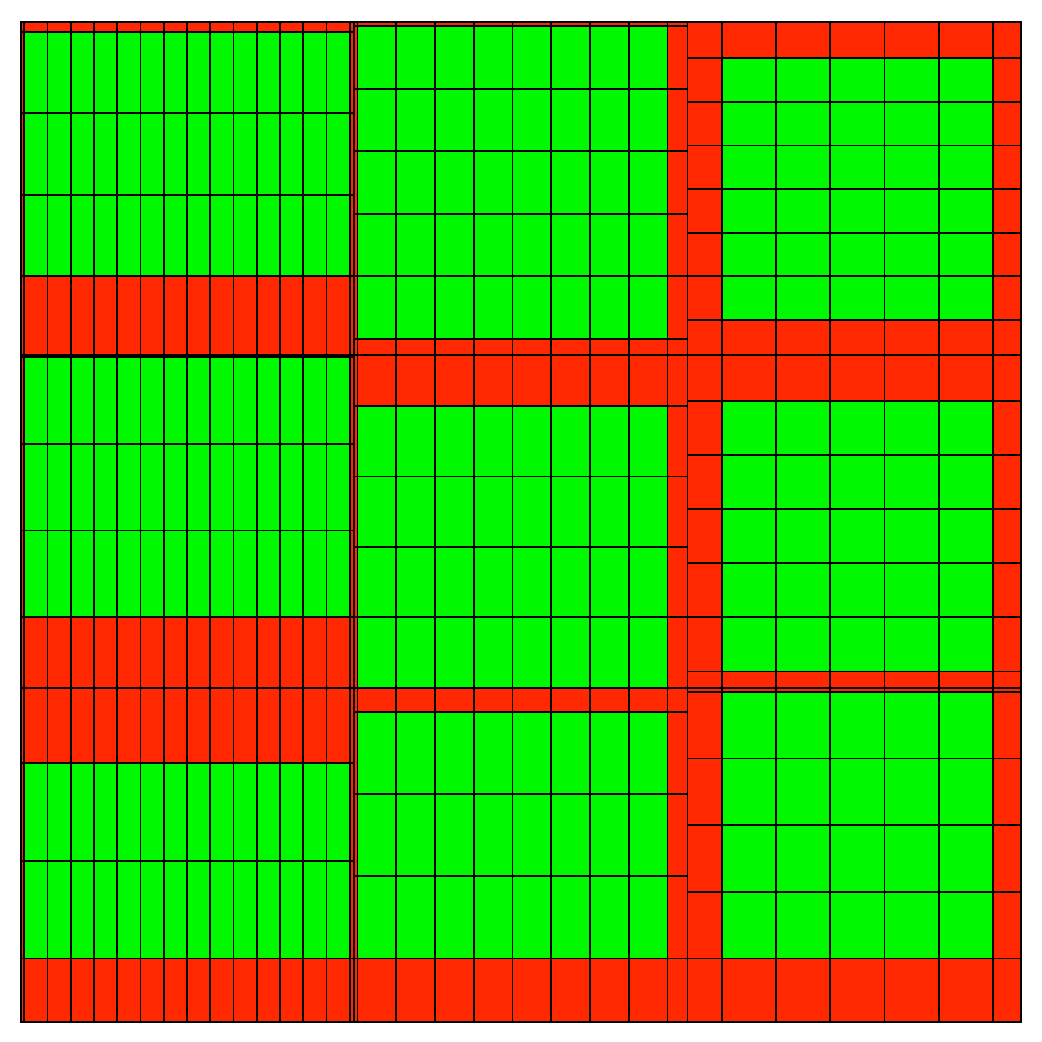}
    \caption{(Left) the initial uniform (coarse) tiling $\cQ_3$ of $R_0$. (Right) the set of blocks $\cP_n^1$ in green and the set of blocks $\cP_n^{2*}$ in red.}
\end{figure}

\end{document}